\documentstyle[12pt]{article}
\begin{document}

\title{Mapping threefolds onto three-quadrics}
\author{Carmen Schuhmann}
\date{}
\maketitle
\newtheorem{lemma}{Lemma}
\newtheorem{remark}{Remark}
\newtheorem{theorem}{Theorem}
\newtheorem{corollary}{Corollary}

\begin{abstract}
We prove that the degree of a nonconstant morphism from a smooth
projective 3-fold $X$ with N\'{e}ron-Severi group ${\bf Z}$
to a smooth 3-dimensio- nal quadric is bounded in terms of numerical
invariants of $X$. In the special case where $X$ is a 3-dimensional
cubic we show that there are no such morphisms. The main tool in the proof is
Miyaoka's bound on the number of double points of a surface.
\end{abstract}

\section {\bf{Introduction}}

During the last years, a lot of progress has been made concerning the
classification of algebraic varieties of dimension larger than 2.
In order to complete the picture, it seems interesting to study the
possible morphisms between these varieties too. Here, only a few
results are known. For instance, Remmert and Van de Ven showed in
\cite{R-V} that the 2-dimensional complex projective space ${\bf P}^2$ cannot
be mapped onto any other smooth complex variety of dimension larger than 0,
and conjectured that the same is true for complex projective spaces of
arbitrary dimension. This conjecture was proven by Lazarsfeld
(see \cite{La}), who used Mori's characterization of ${\bf P}^n$ as the
only smooth, n-dimensional variety with ample tangent bundle.
Later, Paranjape and Srinivas proved that every nonconstant morphism from a
complex homogeneous space, which is not a projective space, to itself
is an isomorphism (see \cite{P-S}). Furthermore, they showed that the complex,
n-dimensional quadric can only be mapped onto ${\bf P}^n$ and itself.
This last result has been generalised to arbitrary characterisic except
2 by Cho and Sato (see \cite{C-S}).

In this paper we work only with varieties and morphisms over
the complex numbers. We will study morphisms from smooth
projective threefolds with N\'{e}ron-Severi group ${\bf Z}$ to smooth
quadric hypersurfaces in ${\bf P}^4$.
The main new tool is the use of Miyaoka's inequality,
bounding the number of double points of a surface with numerically effective
dualizing sheaf in terms of numerical invariants of the surface
(see Theorem \ref{Miyaoka}).
In fact, it will be shown that the
degree of a finite morphism from a smooth
projective threefold $X$ with N\'{e}ron-Severi group ${\bf Z}$ to a smooth
quadric of dimension 3 is bounded in terms of numerical invariants of $X$.
Before giving an exact statement
of this result, we will introduce some definitions for the purpose
of this paper.

Let $X,~Y$ be smooth projective threefolds with N\'{e}ron-Severi group
${\bf Z}$. Denote the ample generator of $NS(X)$ respectively $NS(Y)$ by $H_X$
respectively $H_Y$. The numerical index of $X$ is defined to be the integer $k$
such that the canonical divisor of $X$, $K_X$, is numerically equivalent to
$kH_X$. Let $f$ be a morphism from $X$ to $Y$,
which is not constant. It follows that $f$ is finite (see Lemma \ref{NS}).
The positive integer $d$ which satisfies $f^*H_Y \equiv dH_X$ will be
called the generator degree of $f$. As usual, the number of points, mapping to
a fixed
general point in $Y$ under $f$, will be called the degree of $f$.

The main result of this paper is the following:

\begin{theorem} \label{stelling}
Let $X$ be a smooth projective threefold with N\'{e}ron-Severi group ${\bf Z}$,
$Q$ a
smooth quadric hypersurface in ${\bf P}^4$ and $f$ a nonconstant morphism
from $X$ to $Q$. Then the degree of $f$ is bounded in terms of the
Chernnumbers $c_1^3(X)$ and $c_1^2c_2(X)$ and the numerical index of $X$.
\end{theorem}

This theorem will be proven in Section $2$. We remarked before, that
the main tool, used in the proof, is Miyaoka's inequality.
As Miyaoka's inequality is only valid for surfaces in characteristic 0,
this proof cannot be generalized to other characteristics.

In Section 3, the special case where $X$ is a cubic hypersurface in ${\bf P}^4$
will be treated. We will prove the following result:

\begin{theorem} \label{kubiek}
Let $X$ be a smooth cubic hypersurface in ${\bf P}^4$ and $Q$ a smooth quadric
hypersurface in ${\bf P}^4$. There are no nonconstant morphisms from $X$ to
$Q$.
\end{theorem}

In the proof we will use the results from Section 2, which imply that there
are no morphisms from a cubic to a quadric of generator degree larger than 2.
It will be shown by means of an ad hoc argument, based on
a theorem of Lazarsfeld about morphisms to projective spaces (see Theorem
\ref{Lazarsfeld}), that morphisms of generator degree $2$ cannot occur either.

I thank Prof. Van de Ven for helpful conversations and especially for
suggesting to use \cite{Mi}. I thank Johan de Jong and Endre Szab\'{o} for
stimulating discussions on this subject.

\section{Proof of Theorem 1}

The proof of Theorem \ref{stelling} is based upon the following result of
Miyaoka (see \cite{Mi}):

\begin{theorem} \label{Miyaoka}
Let $S$ be a complex projective surface with only ordinary double points and
numerically effective dualizing sheaf $K_S$. Let $\tilde{S}$ be the minimal
resolution of $S$. Then
$$\#\mbox{\{double points\}} \leq \frac{2}{3}(c_2(\tilde{S}) - \frac{1}{3}
K_S^2).$$
\end{theorem}

Another important ingredient of the proof is the following lemma,
which will be proven later on. Here the tangent hyperplane to a smooth
quadric $Q$ of dimension 3 at a point
$p \in Q$ is denoted by $T_pQ$.

\begin{lemma} \label{lemma}
Let $X$ be a smooth, projective variety of dimension 3 and
$f$ a finite morphism from $X$ to $Q$. Then there is
a dense open subset $U$ in $Q$ such that $f^*(T_pQ \cap Q)$ has no
singularities away from $f^{-1}(p)$ for all $p\in U$.
\end{lemma}

Before starting the proof of Theorem \ref{stelling}, let us state the
following lemma, which has an elementary proof.

\begin{lemma} \label{NS}
Let $X$, $Y$ be smooth projective threefolds and $f:X \rightarrow Y$ a
morphism. Assume that $X$ has N\'{e}ron-Severi group ${\bf Z}$. Then the
following
two statements are equivalent:

$i$) $f$ is nonconstant;

$ii$) $f$ is finite.
\end{lemma}

{\it Proof of Theorem \ref{stelling}:} Suppose $f:X \rightarrow Q$ is
a nonconstant morphism of generator degree d. By Lemma \ref{NS}, $f$ is finite.
The                                 degree of $f$ is equal to
$H_X^3d^3/2$, where $H_X$ is the ample generator of the N\'{e}ron-Severi
group of $X$.
For every point $p$ on $Q$, the hyperplane section $H_p=T_pQ \cap Q$ is a
quadric with one ordinary double point, namely $p$. So, if $p$ is not
contained in the branch divisor of $f$,
then the surface $f^*(H_p)$ contains $H_X^3d^3/2$ ordinary double
points, which map to $p$ under $f$.
{}From Lemma \ref{lemma} it follows that $f^*(H_p)$ has no other
singularities for general $p$.

Now fix a point $p$ such that $f^*(H_p)$ has exactly $H_X^3d^3/2$
ordinary double points and no other singularities and denote $f^*(H_p)$ by
$S_p$. It will be shown that, if $d$ is large enough, Theorem \ref{Miyaoka} can
be applied to $S_p$ and provides an upper bound on the number of ordinary
double points of $S_p$ which is smaller than $H_X^3d^3/2$. In order to
compute this bound, we have to compute $c_2(\tilde{S_p})$ and $K_{S_p}^2$,
where $\tilde{S_p}$ is the minimal resolution of $S_p$.

By Bertini's Theorem (see \cite{Jo}, Th\'{e}or\`{e}me 6.10), there is
a hyperplane section $H$ of $Q$ such that the surface $f^*(H)$ is
nonsingular. Denote this surface by $S$.
The surfaces $S$ and $\tilde{S_p}$ are homeomorphic (see \cite{A}, Theorem 3),
so $c_2(\tilde{S_p})=c_2(S)$ and $c_1^2(\tilde{S_p})=c_1^2(S)$. As $S_p$
has only ordinary double points, it follows that $K_{S_p}^2=K_{S}^2$. Using
$K_X \equiv kH_X$, where k is
the numerical index of $X$, the adjunction formula gives:
$$K_{S} \cong (K_X + S) | _{S} \equiv (k+d)H_X | _{S},$$
so
$$ K_{S}^2 = (k+d)^2H_X^2 | _{S} = (k+d)^2dH_X^3. $$
The second Chernclass of $S$ can easily be computed by means of adjunction:
$$ c_2(S)=dc_2(X)H_X+d^2(d+k)H_X^3.$$
Thus the expression $\frac{2}{3}(c_2(\tilde{S_p}) - \frac{1}{3}
K_{S_p}^2)$, which equals $\frac{2}{3}(c_2(S) - \frac{1}{3}K_S^2)$ by
previous remarks, becomes the following polynomial expression in $d$:
\begin{equation}
\frac{4}{9}H_X^3d^3 + \frac{2}{9}kH_X^3d^2 + \frac{2}{3}(c_2(X)H_X -
\frac{1}{3}k^2H_X^3)d. \label{eq:bound}
\end{equation}
We can apply Theorem \ref{Miyaoka} to $S_p$ if $S_p$ has only double points and
$K_{S_p}$ is nef. As we remarked before, the first condition is certainly
satisfied because of Lemma \ref{lemma}. As for the second one, using
the adjunction formula, it follows that $K_{S_p}$ is linearly equivalent to
$(K_{X}+S_p) |_{S_p}$. So $K_{S_p}$ is nef if and only if $(K_{X}+S_p)
|_{S_p}$ is nef. As $(K_{X}+S_p) |_{S_p} \equiv (k+d)H_X |_{S_p}$,
this is certainly true if $d \geq -k$.
Notice that this condition is empty if $k \geq -1$. If $k=-2$, then the
condition becomes $d \geq2$, which is also an empty condition as there are no
morphisms of degree 1 between a variety of numerical index -2 and the
quadric $Q$, which has numerical index -3. The only smooth threefolds
with numerical index less than -2 are ${\bf P}^3$,
which has numerical index -4 and $Q$ (see \cite{K-O}). So,
the only cases in which
we cannot apply Theorem \ref{Miyaoka} to $S_p$ are when $X$ is ${\bf P}^3$
and $d \leq 3$ or $X$ is $Q$ and $d \leq 2$. In the other cases, it
 tells us that the
number of ordinary double points on $S_p$ is restricted by
the expression (\ref{eq:bound}).
However, we remarked that $S_p$ contains exactly $H_X^3d^3/2$
double points. As the leading term of (\ref{eq:bound}), $4H_X^3d^3/9$,
is smaller than $H_X^3d^3/2$, the expression (\ref{eq:bound})
becomes smaller than $H_X^3d^3/2$ for large $d$. So we obtain a
contradiction if $d$ is larger than the largest positive zero of the polynomial
$H_X^3d^3/2-$(\ref{eq:bound}) (if $X$ is ${\bf P}^3$ resp. $Q$, then we
moreover have to require $d \geq 4$ resp. $d \geq 3$). We conclude
that the generator degree $d$ of $f$ and thus also the degree of $f$
is bounded in terms of the coefficients of this polynomial.
Since $H_X=-c_1(X)/k$, the statement of the theorem follows.\\

{\it Proof of Lemma \ref{lemma}:} Denote the 4-dimensional dual projective
space by ${{\bf P}^4} ^{\vee}$. Let ${Q}^{\vee}$ be the dual variety of $Q$.
It is isomorphic to $Q$ via the following isomorphism:
\[ \begin{array}{rcl}
Q & \longrightarrow & {Q}^{\vee} \\
p & \longmapsto & T_p Q.
\end{array} \]
Given $T_p Q \in {Q}^{\vee}$, denote the surface $f^*(T_p Q \cap Q)
\subset Y$ by $S_p$. We will study the singularities of $S_p$, using
the following criterion:
$$
{}~~~~~~~~~~~~~~~~~~~ x \in S_p \mbox{ is a singularity of } S_p
\Leftrightarrow T_p Q \supset f_* T_x X .~~~~~~~~~~~~~~~~~~~~~~(*)
$$
In order to describe $f_* T_x X$, we introduce some notation. For $i \in \{0,
1,2,3\}$, the set
$$ X_i := \{x \in X \mid f \mbox{ has rank at most } i \mbox{ at }x\} $$
is an algebraic subset of $X$ of dimension $i$. So its image under the
finite map $f$ is an algebraic subset of $Q$ of dimension $i$. Denote
the surface $f(X_2)$ by $B$ (this is just the branch locus of $f$), the
curve $f(X_1)$ by $C$ and the finite set $f(X_0)$ by $R$. Furthermore,
denote the union of all irreducible curves along which $B$ is singular
by $\Gamma$ and the set of isolated singularities of $B$ by $\Sigma$.
Finally, the singular locus of the curve $C$ respectively $\Gamma$ is denoted
by
$Sing(C)$ respectively $Sing(\Gamma)$. Let us now examine when a point $x \in
X$
is a singularity of $S_p$.

If $x \in X_i \backslash X_{i+1}~(i \in \{1,2,3\})$ and $f(x_i)$ is a smooth
point of $f(X_i)$, then $f_*T_xX=T_{f(x)}f(X_i)$. So, according to
criterion $(*)$, $x$ is a singularity of $S_p$ if and only if $T_pQ$
contains $T_{f(x)}X_i$, in other words if and only if $T_pQ$ is tangent
to $f(X_i)$ at the point $f(x)$. Especially, taking $i=3$, it follows that
$x \in X_3 \backslash X_2$ is a singularity of $S_p$ if and only if
$f(x)=p$.

If $x \in X_2 \backslash X_1$ and $f(x) \in \Gamma \backslash Sing(\Gamma)$,
then $f_*T_xX$ contains $T_{f(x)}\Gamma$. So in this case we see from $(*)$
that
if $x$ is a singularity of $S_p$ then $T_pQ$ contains $T_{f(x)}\Gamma$, which
means that $T_pQ$ is tangent to $\Gamma$ at $f(x)$.

Finally, if $x \in X_0$, then $f_*T_xX=0$, so by $(*)$ the point $x$ is
certainly a singularity of $S_p$.

Combining these observations, it follows that $S_p$ has no singularities
away from $f^{-1}(p)$ if the hyperplane $T_pQ$ lies in the intersection
of the following 4 subsets of ${Q}^{\vee}$:

\[
\begin{array}{lll}
U_1 & := & \{H \in {Q}^{\vee} \mid H \cap (R \cup Sing(C) \cup
Sing(\Gamma) \cup \Sigma) = \emptyset \}, \\
U_2 & := & \{H \in {Q}^{\vee} \mid H \mbox{ intersects }
C \mbox{ transversally} \}, \\
U_3 & := & \{H \in {Q}^{\vee} \mid H \mbox{ intersects } \Gamma
 \mbox{ transversally} \}, \\
U_4 & := & \{H \in {Q}^{\vee} \mid H \cap B \backslash (\Gamma \cup
\Sigma) \mbox{ is nonsingular} \}.
\end{array}
\]

Thus, in order to prove the lemma, it is sufficient to show that the sets
$U_i$ ($i \in \{1,2,3,4\}$) are Zariski-open in ${Q}^{\vee}$ and nonempty.

As $R \cup Sing(C) \cup
Sing(\Gamma) \cup \Sigma$ is a finite set, this is certainly true for $U_1$.

As for $U_2$, notice that in order to show that a general
hyperplane $H \in {Q}^{\vee}$ intersects $C$
transversally, it is sufficient to show that this is true for every irreducible
component of $C$. So we can without loss of generality assume that $C$ is
irreducible.

Let $Z$ denote the following closed subscheme of
$Q \times {Q}^{\vee}$:
$$ Q \times {Q}^{\vee} \supset Z := \{(x,H) \in Q \times
{Q}^{\vee} \mid x \in H \}. $$
Consider the scheme-theoretic intersection
$(C \times {Q}^{\vee}) \cap Z$ and the restriction of
the projections from $Q \times {Q}^{\vee}$ onto $Q$ respectively ${Q}^{\vee}$
to this subscheme of $Q \times {Q}^{\vee}$:

\[
\begin{array}{l}
q: (C \times {Q}^{\vee}) \cap Z \longrightarrow C,\\
r: (C \times {Q}^{\vee}) \cap Z \longrightarrow {Q}^{\vee}.
\end{array}
\]

The fiber of $r$ over $H \in {Q}^{\vee}$ is the scheme-theoretic
intersection of $C$ and $H$. So in order to show that $U_2$ is open and dense
in ${Q}^{\vee}$,
we have to prove that the general fiber of $r$ is reduced.
Notice that all fibers of $q$ are singular quadric surfaces, so they
have constant Hilbertpolynomial. It follows that
$q$ is flat (see \cite{Ha}, Chapter III, Theorem 9.9). As all fibers of $q$ are
reduced and $q$ itself is flat, we conclude that
$U := q^{-1}(C \setminus Sing(C))$ is reduced.
Now $U$ is a dense open subscheme of $(C \times {Q}^{\vee}) \cap Z$, so
the general fiber of $r$, which is finite, is contained in $U$.
Restricting to those fibers which are contained in the smooth part
of $U$ and have empty intersection with the ramification locus
of $r$, we see that the general fiber of $r$ is reduced.
It follows that $U_2$ is Zariski-open in ${Q}^{\vee}$ and nonempty.

Replacing $C$ by $\Gamma$ and
repeating the above reasoning proves that $U_3$ is also Zariski-open
in ${Q}^{\vee}$ and nonempty.

In order to prove that $U_4$ is open and dense in ${Q}^{\vee}$, we will show
that, for a general element $H \in
{Q}^{\vee}$, the intersection of $H$ and $B$ is smooth away from the
singular locus $\Gamma \cup \Sigma $ of $B$. As it is sufficient to show that
this holds for every irreducible component of $B$, we may without
loss of generality assume that $B$ is irreducible.

Assume first that $B$ is a hyperplane section of $Q$, corresponding to
an element $H_B$ in ${{\bf P}^4} ^{\vee}$. Then the hyperplane sections of
$B$ correspond to the lines in ${{\bf P}^4} ^{\vee}$ through $H_B$. If
$H_B$ is not contained in ${Q}^{\vee}$, then every line through $H_B$
intersects ${Q}^{\vee}$ in 2 (not necessarily distinct) points. So
every hyperplane section of $B$ can be written as the intersection
of $B$ and an element of ${Q}^{\vee}$. If $H_B$ is an element of
${Q}^{\vee}$, then every line through $H_B$ which is not contained in the
tangent space to ${Q}^{\vee}$ at $H_B$ intersects ${Q}^{\vee}$ in exactly
one point different from $H_B$. So in this case the elements of a dense open
subset of the space of hyperplane sections of $B$ can uniquely be written as
the intersection of $B$ and an element of ${Q}^{\vee}$. In both
cases we conclude from Bertini's Theorem (see \cite{G-H}, page 137),
applied to $B$, that the intersection of $B$ and a general element of
${Q}^{\vee}$ is smooth away from $\Gamma \cup \Sigma $. So
$U_4$ is Zariski-open in ${Q}^{\vee}$ and nonempty.

{}From now on, assume that $B$ is not a hyperplane section of $Q$ and
interprete ${Q}^{\vee}$ as a quadric system contained in the space
${\bf P}(H^0(B,{\cal O}_B(1)))$ of all hyperplane sections of $B$.
In order to prove that the intersection of $H$ and $B$ is smooth away from
$\Gamma \cup \Sigma $ for general $H$ in ${Q}^{\vee}$, we will use the
following special case of Bertini's Theorem (see \cite{G-H},
page 137). Here a one-dimensional linear system is called a pencil.

\begin{lemma} \label{Bertini}
Let $X \subset {\bf P}^N$ be a projective variety and $\Lambda \subset {\bf P}(
H^0(X,{\cal O}_X(1)))$ a pencil. Then the general element of $\Lambda$
is smooth away from the base locus of $\Lambda$ and the singular locus of $X$.
\end{lemma}

We will show that ${Q}^{\vee}$ contains enough pencils to globalise
Lemma \ref{Bertini}, which holds for any of these pencils, to a
Bertini type theorem for the whole space ${Q}^{\vee}$.

For every element $H$ of ${Q}^{\vee}$ there is a one-dimensional family
of pencils in ${Q}^{\vee}$ containing $H$, parametrised by a
plane quadric curve. By Lemma \ref{Bertini}, the general element of a
pencil is smooth away from the base locus of the pencil and
$\Gamma \cup \Sigma $. We will call such an element general for that pencil.
Denote the base locus in $B$ of a pencil $\Lambda$ in ${Q}^{\vee}$
by $B_{\Lambda}$.
Denote the isomorphism from $Q$ to ${Q}^{\vee}$, mapping $p \in Q$ to the
hyperplane $T_p Q$, by $T$. An easy computation shows that
\begin{equation}
B_{\Lambda} = T^{-1}(\Lambda) \cap B. \label{eq:base locus}
\end{equation}
By a dimension argument, there is a Zariski-open subset $V$ of ${Q}^{\vee}$
such that every element $H$ of $V$ is general for almost all
pencils containing $H$. We will show that $V \backslash (V \cap T(B))$ is
contained in $U_4$.

Let $T_p Q$ be an element of $V$. By (\ref{eq:base locus}) the intersection
$B_{\Lambda} \cap B_{\Lambda'}$ of the base loci of any 2 pencils $\Lambda$
and $\Lambda'$ in ${Q}^{\vee}$ is empty if $p \not \in B$ and equals $p$
if $p \in B$. As $T_p Q$ is an element of $V$, it follows that the hyperplane
section $T_p Q \cap B$ is smooth away from $\Gamma \cup \Sigma $ if
$p \not \in B$. So the nonempty Zariski-open subset $V \backslash
(V \cap T(B))$ of ${Q}^{\vee}$ is contained in $U_4$. We conclude that $U_4$
is open and dense in ${Q}^{\vee}$.

\begin{remark} {\rm Notice that, if $X$ is ${\bf P}^3$ respectively a smooth
threedimensional quadric $Q$, then the statement of Theorem \ref{stelling}
follows also from Theorem \ref{Lazarsfeld} respectively the following theorem
(see \cite{P-S}):

\begin{theorem} \label{Pasri}
Let $Y$ be a smooth quadric hypersurface of dimension at least 3, $X$ a smooth
projective variety and $f:Y \rightarrow X$ a surjective
morphism. Then either $f$ is an isomorphism or $X$ is isomorphic to a
projective space.
\end{theorem}

Conversely, we conclude from Theorem \ref{stelling} that the degree of a
nonconstant morphism from $Q$ to itself is bounded. So every such morphism
must have degree 1, as otherwise we could produce nonconstant selfmaps of $Q$
of arbitrarily high degree by composition.

More generally, let $X$ be a smooth threefold with N\'{e}ron-Severi
group ${\bf Z}$. From
Theorem \ref{stelling} it follows that, if $X$ does allow some nonconstant
morphism to the quadric (for example if $X$ is a Fermat hypersurface of
even degree in ${\bf P}^4$), then every nonconstant morphism from $X$ to itself
has degree 1, so is an isomorphism.}
\end{remark}

\begin{remark}{\rm Replacing $Q$ by another smooth threefold with
N\'{e}ron-Severi group
${\bf Z}$, I did not manage to prove an analogue of Theorem \ref{stelling}. For
instance, let $Y$ be a smooth complete intersection of $N-3$ hyperplanes of
degrees $m_1, \dots,m_{N-3}$ in ${\bf P}^N$. Let $X$ be as in Theorem
\ref{stelling}
and $f:X \rightarrow Y$ a finite morphism of generator degree $d$. Then the
degree of $f$ equals $H_X^3d^3/(\prod_{i=1}^{N-3} m_i)$.
So, if $H$ is a hyperplane section of $Y$ with $n$ ordinary double points,
none of which lie in the branch locus of $f$, then $f^*(H)$ has
$nH_X^3d^3/(\prod_{i=1}^{N-3} m_i)$ ordinary double points. As the
leading term of the Miyaoka bound for $f^*(H)$ in Theorem \ref{Miyaoka}
is equal to $\frac{4}{9}H_X^3$, it is clear that the idea of the proof of
Theorem \ref{stelling} can only work in this case if
\begin{equation}
n \geq \frac{4}{9}\prod_{i=1}^{N-3} m_i. \label{eq:ci}
\end{equation}
The point is that, in order to apply Theorem \ref{Miyaoka}, it still has
to be checked that for at least one such hyperplane
section $H$ the surface $f^*(H)$ has only ordinary double
points. If $Y$ is a cubic in
${\bf P}^4$ or the intersection of two quadrics in ${\bf P}^5$, it follows from
(\ref{eq:ci}) that we have to study hyperplane sections of $Y$ with at
least 2 ordinary double points. In both cases I did not succeed in  proving
an analogue of Lemma \ref{lemma}. If $Y$ is a cubic, the system of
hyperplane sections with 2 ordinary double points has only dimension 2 and
I couldn't apply any Bertini type argument; if $Y$ is an intersection of
2 quadrics, this system has dimension 3, but I could not get any grip
on it.

Another possibility is to try to prove something weaker than Lemma \ref{lemma}.
If one can prove that $f^*(H)$ has only mild singularities apart from the
ordinary double points lying over the double points of $H$, then one may try to
apply Theorem \ref{Miyaoka} after blowing up these singularities.
Of course this is only possible if one knows these singularities well.

If one can prove that $f^*(H)$ contains only quotient singularities, then
one may try to apply a more general version of Theorem \ref{Miyaoka}, valid
for surfaces with only quotient singularities (see \cite{Mi}).

However, these approaches seem quite hard and I did not try them
seriously up to now.}
\end{remark}

\section{Proof of Theorem 2}

Look at the special case where $X$ is a smooth hypersurface of degree
$m$ in ${\bf P}^4$ and $f:X \rightarrow Q$ a nonconstant morphism of
generator degree $d$.
As $X$ has N\'{e}ron-Severi group ${\bf Z}$, we can apply Theorem
\ref{stelling}
which gives that
the generator degree of $f$ is bounded. To estimate this bound,
consider expression (\ref{eq:bound}). As

\[ \begin{array}{lll}
H_X^3 & = & m, \\
H_Xc_2(X) & = & m^3 - 5m^2 + 10m
\end{array} \]

and the index of $X$ is equal to $-5+m$, this expression becomes:
\begin{equation}
\frac{4}{9}md^3+(\frac{2}{9}m^2-\frac{10}{9}m)d^2+
(\frac{4}{9}m^3-\frac{10}{9}m^2+\frac{10}{9}m)d. \label{eq:boun}
\end{equation}
The degree of $f$ is equal to $md^3/2$. Thus, in this
case, the upper bound on the generator degree of $f$ we obtained in the proof
of
Theorem \ref{stelling} is given by the maximal positive integer $d$
for which (\ref{eq:boun}) $ \geq md^3/2$. Denote this
integer by $d_m$. A calculation shows that, for $m \gg 0$, the bound $d_m$
we obtain in this way on the generator degree of $f$ grows approximately
linearly with $m$: $d_m \sim (2+2\surd 3)m+constant$.

However, the generator degree of $f$ is also bounded from below.
Choose coordinates $(x_0: \dots :x_4)$ on the projective space ${\bf P}^4$
containing $X$ and coordinates $(y_0: \dots :y_4)$ on the projective
space ${\bf P}^4$ containing $Q$ such that $Q$ is given by the equation
$\sum_{i=0}^{4} y_i^2 =0$. Then $f$ is given by

\[ \begin{array}{rcl}
f:X & \longrightarrow & Q, \\
(x_0: \dots :x_4) & \longmapsto & ({\phi}_0 (x_0: \dots :x_4): \dots :
{\phi}_4 (x_0: \dots :x_4))
\end{array} \]

where the ${\phi}_i$ are homogeneous polynomials of degree $d$, defined on
$X$. As the natural map $H^0({\bf P}^4, {\cal O}_{{\bf P}^4}(d))
\rightarrow H^0(X, {\cal O}_X(d))$ is surjective,
these polynomials can be extended to polynomials of degree $d$, defined on
${\bf P}^4$ (but not necessarily in a unique way). The
extensions will also be denoted by ${\phi}_i$.
Let $X \subset {\bf P}^4$ be given by the equation $F_X=0$, where
$F_X$ is a homogeneous polynomial of degree $m$. As $\sum_{i=0}^{4}
{\phi}_i^2 =0$ on $X$, it follows that
\begin{equation}
\sum_{i=0}^{4} {\phi}_i^2 = F_XG, \label{eq:vgl}
\end{equation}
where $G$ is a homogeneous polynomial of degree $2d-m$. Thus, the
generator degree of $f$ must be larger than or equal to $m/2$.

In fact, if $m$ is even, then there exist hypersurfaces $X$ of degree $m$ and
morphisms of generator degree $m/2$ from $X$ to $Q$, for instance if $X$
is the Fermat hypersurface of degree $m$ in ${\bf P}^4$. From (\ref{eq:vgl}) it
follows that, if $m$ is even, there is a morphism of generator degree $m/2$
from
$X$ to $Q$ if and only if $F_X$ can be written as the sum of 5 squares
of homogeneous polynomials of degree $m/2$, having no common zeroes
on $X$.

Now consider the case $m=3$, where $X$ is a smooth cubic. In order to prove
Theorem \ref{kubiek}, we will use the following theorem of Lazarsfeld (see
\cite{La}):

\begin{theorem} \label{Lazarsfeld}
Let $X$ be a smooth projective variety of dimension at least 1 and let
$f:{\bf P}^n \rightarrow X$ be a surjective morphism. Then $X \cong {\bf P}^n$.
\end{theorem}

{\it Proof of Theorem \ref{kubiek}}: Let $f:X \rightarrow Q$ be a
morphism of generator degree $d$. A computation shows that the
upper bound $d_3$ on the generator degree of $f$, which was introduced
above, is equal to 3. Morphisms of generator degree 3 cannot occur as
the expression $3d^3/2$ which should be equal to the degree of
such a morphism is not integer for $d=3$. So, all that remains to
be proven is that there are no morphisms of generator degree 2
between cubics and quadrics.

Assume $f$ has generator degree 2. As above, choose coordinates
$(x_0: \dots :x_4)$ and $(y_0: \dots :y_4)$ on ${\bf P}^4$, and let
${\phi}_0, \dots, {\phi}_4$ be homogeneous polynomials of degree $2$,
defining $f$. In this case we get, as in (\ref{eq:vgl}):
\begin{equation}
\sum_{i=0}^{4} {\phi}_i^2 = F_XL, \label{eq:kvgl}
\end{equation}
where $L$ is a homogeneous linear polynomial, defining a hyperplane in ${\bf
P}^4$.
This hyperplane will also be denoted by $L$, for convenience.

We claim that the ${\phi}_i$ do not have any common zeroes on the hyperplane
$L$. As the ${\phi}_i$ do not have any common zeroes on $X$, the claim
follows for points in $X \cap L$. Now let $p$ be a point in $L \backslash
(X \cap L)$. If ${\phi}_i(p)=0$ for all $i \in \{0, \dots,4 \}$, equation
(\ref{eq:kvgl}) implies that

$$ \frac{\partial F_XL}{\partial x_i}(p)=0,\mbox{ for all } i \in
\{0, \dots,4 \}. $$

As $L(p)=0$ and $F_X(p)\not = 0$ by assumption, we get:

$$ \frac{\partial L}{\partial x_i}(p)=0 \mbox{ for all } i \in \{0, \dots,4 \}.
$$

But this is impossible because $L$, being a hyperplane, is nonsingular. This
proves the claim.

Thus, the ${\phi}_i$ define a morphism from $L$ to $Q$. Restricted to the
surface $L \cap X$, this morphism equals $f | _{L \cap X}$, so it is
not constant. This is a contradiction by Theorem \ref{Lazarsfeld}.

\begin{remark}
{\rm Notice that the argument in the proof of Theorem \ref{kubiek} works for
every hypersurface $X$ of degree $m$ in ${\bf P}^4$ with a morphism
$f:X \rightarrow Q$ of generator degree $d$ such that $2d-m=1$. So, if $m$ is
odd, there are no morphisms of generator degree $(m+1)/2$ from $X$ to $Q$
(if $m \equiv 1mod4$, this is also clear from the fact that
the expression $md^3/2$ which should be equal to the degree
of such a morphism is not integer in this case).}
\end{remark}

Rijksuniversiteit Leiden, Mathematisch Instituut, Postbus 9512, 2300 RA Leiden,
Nederland; e-mail: schuhman@wi.leidenuniv.nl

\end{document}